\documentclass{ws-mpla}
\usepackage[super]{cite}
\usepackage{graphicx}
\usepackage{subfigure}
\usepackage{slashed}
\usepackage{color}
\usepackage{amsfonts}
\usepackage{amsmath}
\usepackage{extarrows}
\usepackage{amssymb}
\usepackage{verbatim}
\usepackage{epsfig}
\usepackage{todonotes}
\begin{document}

\markboth{Jing-Hui Huang}
{A Fast Numerical solution of the quark's Dyson-Schwinger equation with Ball-Chiu vertex}

\catchline{}{}{}{}{}

\title{A Fast Numerical solution of the quark's Dyson-Schwinger  equation with Ball-Chiu vertex}

\author{\footnotesize Jing-Hui Huang}
\address{School of Institute of Geophysics and Geomatics, China University of Geosciences(Wuhan), Lumo Road 388\\
 Wuhan, 430074, China\\}
 
\author{\footnotesize Xue-ying Duan}
\address{School of Automation, China University of Geosciences(Wuhan), Lumo Road 388\\
 Wuhan, 430074, China\\}
 
 \author{\footnotesize Xiang-yun Hu}
\address{School of Institute of Geophysics and Geomatics, China University of Geosciences(Wuhan), Lumo Road 388\\
 Wuhan, 430074, China\\}

\author{ H. Chen}
\address{School of Mathematics and Physics, China University of Geosciences(Wuhan), Lumo Road 388\\
Wuhan, 430074, China \\huanchen@cug.edu.cn}

\maketitle

\pub{Received (Day Month Year)}{Revised (Day Month Year)}

\begin{abstract}
In this paper, we present two feasible and efficient methods to numerically solve the quark's Dyson-Schwinger(qDSE), the qDSE is mathematical systems of nonlinear integral equations of the second kind with high degrees of freedom. It is difficult to analytically solve the qDSE due to its non-linearity and the singularity. Normally we discrete the singular integral equation by Gauss Legendre integral integration formula, then the approximate solutions of integral equation are obtained by iterative method. The main difficulty in the progress is the unknown function, which is the quark's propagator at vacuum and at finite chemical potential, occurs inside and outside the integral sign. Because of the singularity, the unknown function inside the integral sign need to be interpolate with high precision. Normally traditional numerical examples show the interpolation will cost a lot of CPU time. In this case, we provide two effective and efficient methods to optimize the numerical calculation, one is we put forward a modified interpolation method to replace the traditional method. Besides, the technique of OpenMP and automatic parallelization in GCC is another method which has widely used in modern scientific computation. Finally, we compare CPU time with different algorithm and our numerical results show the efficiency of the proposed methods.

\end{abstract}

\ccode{PACS Nos.: 25.75.Nq, }

\section{Introduction}	

The Dyson-Schwinger Equations (DSEs) of QCD are the equations of motion in the continuum quantum field theory.
They are a set of coupled equations about quark, gluon, and ghost propagators and vertex functions.
They provide a successful description of hadrons at vacuum (Nambu phase) (see, e.g., Refs.~\cite{Maris2003,CLR,GCL2017})
and phase transitions in hot/dense medium (see, e.g., Refs.~\cite{QCDPT-DSE11,QCDPT-DSE12-5,QCDPT-DSE23}). However, the previous works
are phenomenological models with effective gluon and vertex input, because the computing power is far beyond the reach of the full numerically solving Dyson-Schwinger Equations with quarks, gluons and vertexes. In this case, it is imperative to come up with more efficient methods to improve computational efficiency.

The qDSE is reduced to mathematically solving the systems of nonlinear integral equations of the second kind with high degrees of freedom. It is difficult to analytically solve the qDSE due to its non-linearity and singularity. Fredholm integral equations occurs in many branches of scientific fields\cite{1973Numerical,0Introduction,linandonlin}, such as microscopy, radio astronomy, electron emission and optical fiber evaluation. The systems of nonlinear integral equations of the second kind reads:
\begin{eqnarray}
u(x)=f_{1}(x)+\int [ K_{1}(x,t) F_{11}(u(t),v(t))+ \overline{K}_{1}(x,t) F_{12}(u(t),v(t))]dt,  \\ 
v(x)=f_{2}(x)+\int [K_{2}(x,t) F_{21}(u(t),v(t))+ \overline{K}_{2}(x,t) F_{21}(u(t),v(t))]dt. \notag
\label{jifen1}
\end{eqnarray}
where $F_{ij}(u(x),v(x))(i=1,2)$ is a nonlinear function, the kernels $K_1$ and $K_2$ are the nonsingular kernels, while the kernels $\overline{K_1}$ and $\overline{K_2}$ are the singular kernels given by 
\begin{eqnarray}
\overline{K}_{1,2}=\frac{1}{(x-t)^{\alpha}}
\label{jifen3}
\end{eqnarray}
The integral will be singular when the eq(\ref{jifen3}) becomes zero at one or more points at the range of integration.It is difficult to analytically solve the system(Eq.\ref{jifen1}) due to its non-linearity and the singularity. Normally we discrete the singular integral equation by Gauss Legendre integral integration formula, then the approximate solutions of integral equation are obtained by iterative method. Specifically the successive approximations method\cite{linandonlin} introduces the recurrence relation:
\begin{eqnarray}
u_{0}(x)=any\ selective\ real \ valued \ function  \nonumber \\
v_{0}(x)=any\ selective\ real \ valued \ function.
\label{jifen14}
\end{eqnarray}
\begin{eqnarray}
u_{n+1}(x)=f_{1}(x)+\int [K_{1}(x,t) F_{11}(u_{n}(t),v_{n}(t))+ \overline{K}_{1}(x,t) F_{12}(u_{n}(t),v_{n}(t))]dt,  \nonumber  \\
v_{n+1}(x)=f_{2}(x)+\int [K_{2}(x,t) F_{21}(u_{n}(t),v_{n}(t))+ \overline{K}_{2}(x,t) F_{21}(u_{n}(t),v_{n}(t))]dt.
\label{jifen12}
\end{eqnarray}
Consequently, the solution
\begin{eqnarray}
u(x)= \mathop{lim }\limits_{n \to \infty} u_{n}(x),  \nonumber  \\
v(x)= \mathop{lim }\limits_{n \to \infty} v_{n}(x).
\label{jifen122}
\end{eqnarray}

The main difficulty in the progress is the unknown function $u(x)$ and $v(x)$, occurs inside and outside the integral sign. Because of the singularity, the unknown function inside the integral sign need to be interpolate with high precision. In the traditional interpolation progress\cite{Peckover_1971}, such as Newton interpolation, spline interpolation, the step to find the location where the unknown interpolation point on the discrete data will cost a lot of CPU times.

Therefore, we intend to provide two effective and efficient methods to optimize the numerical calculation, one is we put forward a modified interpolation method to replace the traditional method. Besides, the technique of Open Multiprocessing (OpenMP) is another method which has widely used in modern scientific computation. OpenMP has been very successful in exploiting structured parallelism in applications\cite{1998OpenMP,2001Parallel,2009The}. Particularly article\cite{openmpgcc} describes the design  of the OpenMP specification v2.5 in GCC. The implementation supports all the languages specified in the standard (C, C++ and Fortran), and it is generally available on any platform that supports POSIX threads. In this case, we use the OpenMP and automatic parallelization in GCC to accelerate our C program.  

The paper is organized as follows. In section II, the truncation scheme of DSE for quark propagator at vacuum is given. In section III, we briefly describe our modified algorithm to sole systems of nonlinear integral equations.
Then, the numerical results are given in section IV. Finally, we summarize our work and give a brief remark in section V.

\section{Dyson-Schwinger Equation for the quark propagator }
\label{SecII}
Our calculation is based on the quark propagator at vacuum,
which satisfies the Dyson-Schwinger equation
\begin{eqnarray}
S(p ; \mu)^{-1}&= &Z_2 (i\gamma\cdot \tilde{p}+m_{q}) 
 +  Z_1 g^2(\mu)\int \frac{d^4q}{(2\pi)^4}
 \nonumber \\ & &\times\!D_{\rho\sigma}(k;\mu) \frac{\lambda^a}{2}\gamma_\rho S(q;\mu)
\Gamma^a_\sigma(q,p;\mu) \, ,
\label{gendse}
\end{eqnarray}
where $k=p-q$, $D_{\rho\sigma}(k ;\mu)$ is the full gluon propagator, $\Gamma^a_\sigma(q,p;\mu)$ is the dressed quark-gluon vertex, $Z_{1}$ is the renormalization constant for the quark-gluon vertex, and $Z_{2} $ is the quark wave-function normalization constant. The general structure of the quark propagator at vacuum can be written as
\begin{equation}
S^{-1}(p) = i\gamma\cdot pA(p^2) + B(p^2) \, .
\end{equation}

We solve Eq.~(\ref{gendse}) with models of the gluon propagator and the quark-gluon vertex, which describe meson properties at vacuum well in the symmetry-preserving Dyson-Schwinger equation and Bethe-Salpeter equation (BSE) schemes (see, e.g., Refs.~\cite{Fischer2009,Chang2009}).
At vacuum they are usually taken as
\begin{equation}
 { Z_{1}g^2 D_{\rho \sigma}(k) \Gamma_\sigma^a(q,p)}  =
{\cal G}(k^2)D^0_{\rho
\sigma}(k)\frac{\lambda^a}{2}\Gamma_{\sigma}(p,q) \, ,
\label{KernelAnsatz}
\end{equation}
where $D^{0}_{\rho \sigma}(k)=\frac{1}{k^2}\Big[\delta_{\rho\sigma}-\frac{k_\rho k_\sigma}{k^2} \Big]$ is the Landau-gauge free gauge-boson propagator, ${\cal G}(k^2)$ is a model effective interaction, and $\Gamma_{\sigma}(q,p)$ is the effective quark-gluon vertex. 

In our work,  we use the the Ball-Chiu (BC) vertex ans\"{a}tz which satisfies the nonperturbative Ward-Takahashi identity~\cite{BC19811,Chen2008}. The BC vertex at vacuum was given in Ref.\cite{BC19811}
\begin{eqnarray}
\Gamma_\sigma^{BC}(q,p;\mu) = && \lambda_{1} \gamma_{\mu} \nonumber \\
&&+\lambda_{2} (p+q)_{\mu}  \nonumber \\
&&+\lambda_{3} (p+q)_{\mu}(p\cdot\gamma+q\cdot\gamma) \nonumber \\
&&+\lambda_{4} (p+q)_{\nu} \delta_{\mu \nu}. \label{bcvtxmu}
\end{eqnarray}
where $k=q-p$,$t=q+p$,$\lambda_{i}$(i=1,2,3,4) is the function of the scaler function $A(p^{2})$,$A(q^{2})$,$B(q^{2})$ and $B(p^{2})$ of the quark's propagator:

\begin{align}\label{2}
 &\lambda_{1}=\frac{1}{2}\frac{A(p^{2}+A(q^{2})}{2},\nonumber \\
 &\lambda_{2}=-i\frac{B(p^{2})-B(q^{2})}{p^{2}-q^{2}},\nonumber \\
 &\lambda_{3}=\frac{1}{2}\frac{A(p^{2})-A(q^{2})}{p^{2}-q^{2}},\nonumber \\
 &\lambda_{4}=0 .
\end{align}

For the model effective interaction, we employ two infrared-dominant models, denoted as the “HF" model,
which only express the long-range behavior of the renormalization-group-improved 
the Qin-Chang(QC) model~\cite{QC2011}. 
Though the ultraviolet parts of the models ensure the correct perturbative behavior,
they are not essential in describing nonperturbative physics, e.g. spectrum of light mesons ~\cite{Eichmann2016}.
Herein we neglect them since our work is based on the competition between the chiral condensate and quark number density in the infrared region ~\cite{Chen2009}.
The "HF" model is expressed as:
\begin{equation}
\label{IRGsQC} 
{\cal G}^{QC}(k^{2}) = \frac{8\pi^2}{\omega^4} \, D\, {\rm e}^{-k^{2}/\omega^2} \, .
\end{equation}

Hence, the regularization mass scale can be removed to infinity and the renormalization constants can be set to 1.

With given the effective quark-gluon vertex and the model effective interaction in the previous section, we could obtain  the systems of nonlinear integral equations of the scaler function $A(p^{2})$ and $B(p^{2})$:
\begin{align}\label{dsesys}
 &A(p^{2})=z_{1}+\int \frac{dq^4}{(2 \pi)^{4}} \frac{{\cal G}^{QC}(k^{2})}
 { k^{2}p^{2} (p^{2}A^{2}(p^2)+B^{2}(p^{2}))}(I_{A1}+I_{A2}+I_{A3})\nonumber \\
 &B(p^{2})=m_{0}z_{1}+\int \frac{dq^4}{(2 \pi)^{4}} \frac{{\cal G}^{QC}(k^{2})}
 { k^{2} (p^{2}A^{2}(p^2)+B^{2}(p^{2}))}(I_{B1}+I_{B2}+I_{B3})
\end{align}
with
\begin{align}
\label{dsesys2}
 & I_{A1}=-\frac{B(q^{2})}{p^2}\frac{k^{2}(p\cdot t)-(k \cdot p)(k \cdot t)}{k^{2} p^{2}} \frac{B(q^{2})-B(p^2)}{q^{2}-p^{2}}\nonumber \\
 & I_{A2}=-\frac{A(q^{2})}{2p^2} \frac{p^{2}(q\cdot t)k^{2}+q^{2}(p\cdot t)k^{2}-q^{2}(p \cdot k)(t \cdot k)-p^{2}(q\cdot k)(t \cdot k)}{k^{2}}  \frac{A(q^{2})-A(p^2)}{q^{2}-p^{2}}\nonumber \\
 & I_{A3}=\frac{A(q^{2})}{p^2} \frac{k^{2}(p \cdot q)+2(k \cdot q)(k \cdot p)}{k^{2}} \frac{A(q^{2})+A(p^2)}{2}\nonumber \\
 & I_{B1}=-A(q^{2}) \frac{(q \cdot t) k^{2} -(q \cdot k)(t \cdot k)}{k^{2}}\frac{B(q^{2})-B(p^2)}{q^{2}-p^{2}} \nonumber \\
 & I_{B2}=3B(q^{2}) \frac{A(q^{2})+A(p^2)}{2} \nonumber \\
 & I_{B3}=B(q^{2}) \frac{k^{2}t^{2}-(k\cdot t)^{2}}{2k^{2}} \frac{A(q^{2})-A(p^2)}{q^{2}-p^{2}}
\end{align}
we could find the singular kernel $1/(q^{2}-p^{2})$ in the term $I_{A1}$,$I_{A2}$,$I_{B1}$ and $I_{B3}$. Normally we discrete the singular integral equation by Gauss Legendre integral integration formula, then the approximate solutions of integral equation are obtained by iterative method with a proper initial function. In our numerical calculation, initial functions are obtained with:
\begin{eqnarray}
A_{0}(p)=1.0,\nonumber \\
B_{0}(p)=1.0.
\label{jifenab1}
\end{eqnarray}
and the solution is determined by the accuracy $\xi$ of convergence at every momentum points:
\begin{eqnarray}
A(p)=\mathop{lim }\limits_{n \to \infty} A_{n}(p) \ \ \longrightarrow   \ \ abs(A_{n+1}(p)-A_{n}(p))< \xi, \nonumber\\
B(p)=\mathop{lim }\limits_{n \to \infty} B_{n}(p) \ \ \longrightarrow  \ \  abs(B_{n+1}(p)-B_{n}(p))< \xi \nonumber.
\label{jifenab2}
\end{eqnarray}

\section{modified interpolation algorithm and automatic parallelization in GCC}

The main difficulty in the numerical solution is the unknown function $A(q^2),B(q^2)$, occurs inside and outside the integral sign and the value of momentum $p$, $q$ always take different value to eliminate the singularity.  Therefore, the unknown function inside the integral sign need to be interpolate with high precision. In addition, higher  calculation accuracy requires more calculation time, and on this base we put forward two effective methods to accelerate our numerical  efficiency.
\subsection{modified interpolation algorithm}
Firstly, we briefly introduce the traditional interpolation and the Computational efficiency.
The traditional interpolation\cite{Peckover_1971}, such as Newton interpolation, spline interpolation, need to find the location of interpolation value point $q_{i}$ on the discrete data p[i](i=1,2,3...N), where the N is the number of the discrete points of the outside momentum.
\begin{gather*}
p[i-1]<q<p[i] \nonumber  \\
 \downarrow where \ ?  \nonumber  \\
p[1]<p[2]< \cdots <p[i-1]<p[i]<p[i+1]< \cdots <p[N]
\label{jifen16}
\end{gather*}
Normally this step (defined "finding step") will cost a lot of CPU time due to logical operations. Then unknown function $F(q^2)$= $A(q^2),B(q^2)$ can be obtained directly:
\begin{equation}
F(q)=\left \{ 
\begin{aligned}
F(p[i])+(q-p[i]) \times \frac{F(p[i+1])-F(p[i])}{p[i+1]-p[i]} \ \ i<N,\\ \\
F(p[i]) \ \ i=N,
\end{aligned}
\right.
\label{chazhiresult0}
\end{equation}

Therefore, we put forward modified interpolation method to replace the traditional one. In particular, the momentum q which defined as integral inner momentum need to be determined beforehand and saved as array q[i](i=1,2,3...N), array p[i] and array q[i] has the relationship:
\begin{equation}
p[1]<q[1]<p[2]< \cdots <p[i-1]<q[i-1]<p[i]< \cdots q[N-1]<p[N]<q[N]
\label{pqrelationship}
\end{equation}
In our modified interpolation method, we use Eq.(\ref{pqrelationship}) to replace the traditional "finding step", then  a lot of CPU time computing logical operations will be saved and the unknown function $F(q^2)$= $A(q^2),B(q^2)$ can be obtained directly:
\begin{equation}
F(q[i])=\left \{ 
\begin{aligned}
F(p[i])+(q[i]-p[i]) \times \frac{F(p[i+1])-F(p[i])}{p[i+1]-p[i]} \ \ i<N,\\ \\
F(p[i]) \ \ i=N,
\end{aligned}
\right.
\label{chazhiresult}
\end{equation}

\subsection{automatic parallelization in GCC}
In the every progress of integration(Eq.\ref{jifen12}), function $F(p[i])$ with outside momentum is certain, while function $F(q[i])$ with inner momentum need to be estimated by interpolation. 

In this work, we use the OpenMP and automatic parallelization in GCC to accelerate our C program by following the standard method from the article\cite{openmpgcc}. Firstly, the pseudo code of a normal "for" loop is showed in the Listing 1:
\begin{center} 
\fcolorbox{black}{gray!10}{\parbox{.9\linewidth}{
void main() \\
\{ \\
  $\quad$ double Ap[N],Bp[N]; \\
  !!... some initialization code for computing integral; \\
  for(int i=0;i<N:i++) \\
  \{ \\
  Ap[i]=functionA(p[i]);\\
  Bp[i]=functionB(p[i]);\\
  \}\\
\}
}}\\
Listing 1. Sequential processing in a for loop 
\end{center}

where functionA(p[i]) and functionB(p[i]) are subroutine computing the right part of Eq.(\ref{dsesys}). Clearly, we can potentially split the "for" loop and run into multiple cores, calculating any Ap[i] is dependency free from other elements of the Ap array, and it is the same to Bp[i]. Then, Listing 2. (pseudo code) shows how OpenMP help us do so. 

\begin{center} 
\fcolorbox{black}{gray!10}{\parbox{.9\linewidth}{
void main() \\
\{ \\
  $\quad$ double Ap[N],Bp[N]; \\
  !!... some initialization code for computing integral; \\
  \# pragma omp parallel for\\
  for(int i=0;i<N:i++) \\
  \{ \\
  Ap[i]=functionA(p[i]);\\
  Bp[i]=functionB(p[i]);\\
  \}\\
\}
}}\\
Listing 2. parallel processing in a for loop with the parallel for pragma
\end{center}

The parallel for pragma helps to split the "for" loop workload across multiple threads, with each thread
potentially running on different cores, thus reducing the total computation time significantly. The exact numerical result will be shown in the next section.

\section{numerical results}
To carry out the numerical calculations, we need the parameters $D$ and $\omega$ in the effective interaction. Usually the parameters are determined by fitting meson properties with the BSE approach. 
The parameters $D=0.550$ and $\omega=0.678$ are taken from  Ref~\cite{PhysRevD.101.054007}. With the above-determined parameters and the ans\"{a}tz, we solve the Dyson-Schwinger equation of the quark propagator with four different algorithms :
algorithm.1) traditional interpolation with sequential processing\\
algorithm.2) modified interpolation with  sequential processing\\
algorithm.3) traditional interpolation with parallel processing \\
algorithm.4) modified interpolation with parallel processing\\
and for each algorithm, there are three parameters to controlling the computational accuracy and CPU time: the length of array p[i](i=1,2..N), the number M of the nodes in Gauss Legendre integral integration formula, and the accuracy $\xi$ of convergence. the Table~\ref{tab:table1} shows our numerical results with various algorithm. 

\begin {table}[h]
\label{tab:table1} 
\tbl{
CPU time With different algorithm and parameters. All codes are run on the computer with Inter(R) Core(TM) i7- @1.00GHz on GCC version 9.3.0 (Ubuntu 9.3.0-10)}
{\begin{tabular}{c@{\hspace{0.7cm}}c@{\hspace{0.7cm}}c@{\hspace{0.7cm}}c@{\hspace{0.7cm}}c@{\hspace{0.7cm}}c@{\hspace{0.7cm}}c@{\hspace{0.7cm}}c@{}}
\toprule
N &$M$ & $\xi$ & algorithm  & \% CPU   & iteration times  & CPU time(s)  \\
\colrule
150  &100&  0.005 &   algorithm.1   & 100.0 &  11 &     65 \\
150  &100&  0.005 &   algorithm.2   & 100.0 &  11 &      6 \\
150  &100&  0.005 &   algorithm.3   & 782.0 &  11 &      9 \\
150  &100&  0.005 &   algorithm.4   & 792.0 &  11 &      1 \\
\botrule
\end{tabular}}
\end{table}

\begin{figure}[h]
\centering
\includegraphics[scale=0.53]{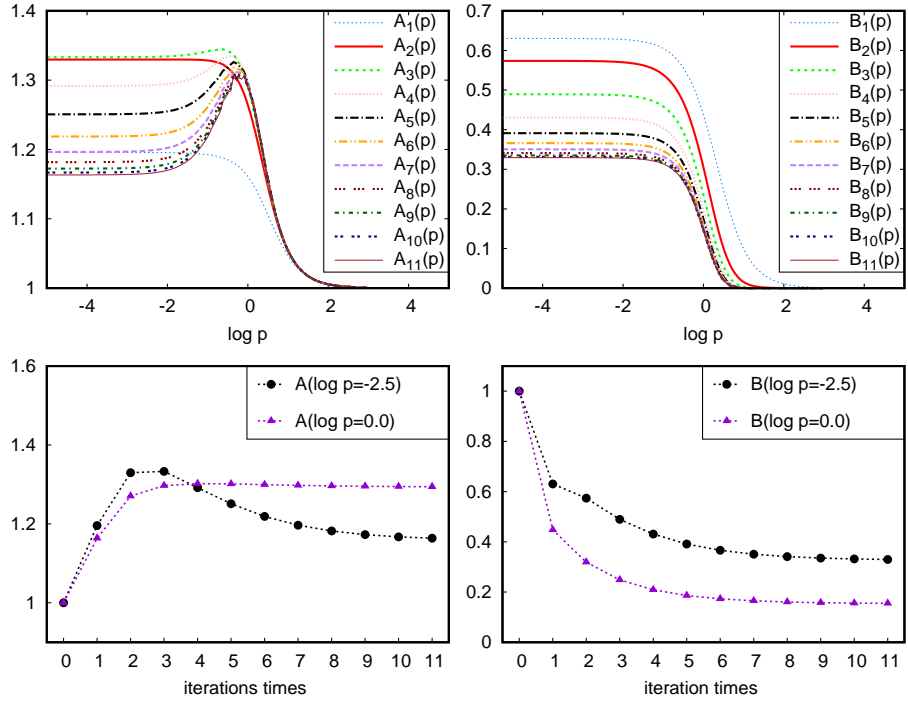}
\caption{ On the up results for the quark's scalar function A(p) and B(p) evolving with iterations. On the up results the quark's scalar function A(log p=-2.5), A(log p=0.0), B(log p=-2.5) and B(log p=0.0) evolving with iterations.}
\label{itresult}
\end{figure}

It is noted that the four algorithms lead to the same solution with same number of iterations. The quark's scalar function A(p) and B(p) in different algorithms consistently evaluates with iterations as shown in Fig.(\ref{itresult}). In a way, the consistently evaluation with different algorithms verify the correctness of our code. Namely, our algorithms just change the efficiency of numerical computation, while the process and the result still remain the same.  
The Table.\ref{tab:table1} shows efficiency of our proposed methods. Both modified interpolation method and the OpenMP and can save a lot of computation time. 

When comparing the result of algorithm.2 and algorithm.1, the modified interpolation method speeds up the traditional one 10 times. When comparing the result of algorithm.3 and algorithm.1, the OpenMP and automatic parallelization speeds up the traditional one about 7 times, which identifies with the computation ability of my computer(Inter i7). When comparing the result of algorithm.4 and algorithm.1, the combine of modified interpolation method and the OpenMP and automatic parallelization speeds up the traditional one about 60 times.

\section{Summary}
In summary, we present two feasible and efficient methods to numerically solve the quark's Dyson-Schwinger(qDSE), which is mathematical systems of nonlinear integral equations of the second kind with high degrees of freedom. The main difficulty in the progress is the unknown function, which is the quark's propagator at vacuum and at finite chemical potential, occurs inside and outside the integral sign. Because of the singularity, the unknown function inside the integral sign need to be interpolate with high precision. In this case, we provide two effective and efficient methods to optimize the numerical calculation, and our numerical results show the efficiency of the proposed methods.

Our interpolation method, which belong to optimization of the code's essence is universal to other interpolation code, while the OpenMP and automatic parallelization will be limited by the hardware configuration of computers. For more complex system, such as the quark's Dyson-Schwinger equation at finite chemical potential or non-zero temperature\cite{Chen2008,QCDPT-DSE12-1,QCDPT-DSE22-1}, the Gauss Legendre integral will be added one more integral degree of freedom, meanwhile two-dimensional interpolation will cost more computation time, hence the combination of the two methods is more efficiently to solve the complex system and the relevant work is undergoing.

\section*{Acknowledgements}

We acknowledge financial support from National Natural Science Foundation of China (Grants No. G1323519204). 

\bibliographystyle{ws-mpla}
\bibliography{main}

\end{document}